\documentclass[12pt,twoside]{article}

\usepackage{amssymb}
\usepackage{amsmath}
\usepackage{bbm}
\usepackage{mathrsfs}
\usepackage{float}
\usepackage{graphicx}
\usepackage{multicol}
\usepackage{pifont}

\newcommand{\br}{ }
\newcommand{\brr}{, }

\makeatletter
\gdef\th@mychange{\normalfont\slshape
   \def\@begintheorem##1##2{\item
        [\hskip\labelsep \theorem@headerfont ##2. ##1  \,--\!--\!--\!--  ]}%
 \def\@opargbegintheorem##1##2##3{%
   \item[\hskip\labelsep \theorem@headerfont ##2. ##1\ {\upshape(}##3{\upshape)}. \,-----  ]}}
\makeatother

\RequirePackage{theorem}
\theoremstyle{mychange}

{\theorembodyfont{\rmfamily}\newtheorem{ttt}{}[section]}
{\theorembodyfont{\rmfamily}\newtheorem{defi}[ttt]{Definition.}}
{\theorembodyfont{\rmfamily}\newtheorem{rem}[ttt]{Remark.}}
{\theorembodyfont{\rmfamily}\newtheorem{rems}[ttt]{Remarks.}}
{\theorembodyfont{\rmfamily}\newtheorem{exa}[ttt]{Example.}}
{\theorembodyfont{\rmfamily}\newtheorem{exas}[ttt]{Examples.}}
{\theorembodyfont{\rmfamily}}
{\theorembodyfont{\rmfamily}\newtheorem{qu}[ttt]{Question.}}
{\theorembodyfont{\rmfamily}\newtheorem{summ}[ttt]{Summary.}}

{\theorembodyfont{\rmfamily}\newtheorem{algso}[ttt]{Algorithms}}

{\theorembodyfont{\itshape}\newtheorem{fac}[ttt]{Fact.}}
{\theorembodyfont{\itshape}\newtheorem{lem}[ttt]{Lemma.}}
{\theorembodyfont{\itshape}}
{\theorembodyfont{\itshape}}

{\theorembodyfont{\itshape}\newtheorem{faco}[ttt]{Fact}}
{\theorembodyfont{\itshape}\newtheorem{lemo}[ttt]{Lemma}}
{\theorembodyfont{\itshape}\newtheorem{propo}[ttt]{Proposition}}
{\theorembodyfont{\itshape}\newtheorem{cono}[ttt]{Conjecture}}

{\theorembodyfont{\rmfamily}}

{\theorembodyfont{\rmfamily}\newtheorem{lui}[ttt]{The method of van Luijk in~detail.}}

\newcommand{\calO}{\mathscr{O}}

\newcommand{\calV}{\mathscr{V}}

\newcommand{\bbC}{{\mathbbm C}}
\newcommand{\bbF}{{\mathbbm F}}

\newcommand{\bbQ}{{\mathbbm Q}}
\newcommand{\bbZ}{{\mathbbm Z}}

\newcommand{\bP}{{\bf P}}
\newcommand{\rk}{\mathop{\text{\rm rk}}}

\newcommand{\Frob}{\mathop{\rm Frob}}
\newcommand{\Pic}{\mathop{\rm Pic}}

\newcommand{\Br}{\mathop{\rm Br}}
\newcommand{\disc}{\mathop{\rm disc}}
\newcommand{\reg}{\text{reg}}
\newcommand{\alg}{\text{alg}}
\newcommand{\et}{\text{\rm \'et}}

\newcommand{\bigtimes}{\mathop{\mbox{\Pisymbol{pzd}{53}}}\limits}

\newcounter{abc}
\newenvironment{abc}{\begin{list}{\rm \alph{abc}) }%
{\usecounter{abc} \leftmargin=0.0pt \labelsep=0.0pt %
\listparindent=0.0pt \labelwidth=0.0pt \parsep=\smallskipamount %
\itemsep=0.0pt \topsep=0.0pt \partopsep=\smallskipamount}}{\end{list}}
\newcounter{iii}
\newenvironment{iii}{\begin{list}{\rm \roman{iii}) }%
{\usecounter{iii} \leftmargin=0.0pt \labelsep=0.0pt %
\listparindent=0.0pt \labelwidth=0.0pt \parsep=\smallskipamount%
 \itemsep=0.0pt \topsep=0.0pt \partopsep=\smallskipamount}}{\end{list}}

\def\rightend#1#2{{%
 \leavevmode\nobreak\hskip .5em plus 1fil
 \penalty600 \hskip 0pt plus -1filll
 \vadjust{}\nobreak\hskip 0pt plus 1filll%
 #1\parfillskip=#2\relax \par}}

\def\eop{\ifmmode\rule[-22pt]{0pt}{1pt}\ifinner\tag*{$\square$}\else\eqno{\square}\fi\else\rightend{$\square$}{0pt}\fi}

\author{
Andreas-Stephan Elsenhans${}^*$ and J\"org Jahnel${}^\ddagger$}

\date{}

\title{Determinantal quartics and the computation of \\the Picard~group}

\begin{document}

\maketitle

\begin{abstract}
We test the methods for computing the Picard group of a
$K3$~surface
in a situation of high~rank. The~examples chosen are resolutions of quartics
in~$\bP^3$
having
$14$
singularities of
type~$A_1$.
Our~computations show that the method of R.\,van Luijk works well
when sufficiently large primes are used.
\end{abstract}

\footnotetext[1]{Mathematisches Institut, Universit\"at Bayreuth, Univ'stra\ss e 30, D-95440 Bayreuth, Germany,\\
{\tt Stephan.Elsenhans@uni-bayreuth.de}, Website:\! {\tt http://www.staff\!.\!uni-bayreuth.de/$\sim$btm216}\smallskip}

\footnotetext[3]{\mbox{Fachbereich 6 Mathematik, Universit\"at Siegen, Walter-Flex-Str.~3, D-57068 Siegen, Germany,} \\
{\tt jahnel@mathematik.uni-siegen.de}, Website: {\tt http://www.uni-math.gwdg.de/jahnel}\smallskip}

\footnotetext[1]{The first author was supported in part by the Deutsche Forschungsgemeinschaft (DFG) through a funded research~project.\smallskip}

\footnotetext[1]{${}^\ddagger$The computer part of this work was executed on the servers of the chair for Computer Algebra at the University of Bayreuth. Both authors are grateful to Prof.~M.~Stoll for the permission to use these machines as well as to the system administrators for their~support.}

\section{Introduction}

\begin{ttt}
The~methods to compute the Picard~rank of a
$K3$~surface~$V$
are limited up to~now. As~shown, for example in~\cite{vL} or~\cite{EJ1}, it is possible to construct a
$K3$~surface
with a prescribed Picard~group. But~when a
$K3$~surface
is given, say, by an equation with rational coefficients, then it is not entirely clear whether its geometric Picard rank may be determined using the methods presently~known.
\end{ttt}

\begin{ttt}
Generally~speaking, it is always possible to give upper and lower~bounds. For~the lower bound, it is necessary to specify divisors explicitly and to verify that their intersection matrix is nondegenerate. This~part is definitely~problematic. It~might happen that a nontrivial divisor is hidden somewhere and very difficult to~find.

On~the other hand, the general strategy for the computation of upper bounds is to use reduction
modulo~$p$.
The~idea to use
characteristic~$p$
methods here is due to R.\,van~Luijk. We~will describe this approach in more detail in~\ref{Luijk}.

Observe,~however, that the Picard rank of a
$K3$~surface
over~$\overline\bbF_{\!p}$
is conjectured to be always~even. In~particular, if
$\rk \Pic(V_{\overline\bbQ})$
is odd then there is no
prime~$p$
such
that~$\rk \Pic(V_{\overline\bbF_{\!p}}) = \rk \Pic(V_{\overline\bbQ})$.
Even~more, the rank
over~$\overline\bbQ$
being even or odd, there is no obvious reason why there should exist a prime
number~$p$
such that
$\rk \Pic(V_{\overline\bbF_{\!p}})$
is at least close
to~$\rk \Pic(V_{\overline\bbQ})$.
\end{ttt}

\begin{ttt}
The~goal of this article is to test van Luijk's method on a randomly chosen sample of
$K3$~surfaces. As~mentioned above, the central point is the existence of good~primes.
Here, being good shall mean that the geometric Picard rank of the reduction
modulo~$p$
does not exceed the Picard rank
over~$\overline\bbQ$
by more than~one.

We~will focus on surfaces of
Picard~rank~$\geq\!15$.
The~reason tor this is a practical~one. For~surfaces of small Picard~rank,
one is forced to work with very small primes such as
$2$
or~$3$
as, otherwise, the computations run out of~time. This~would make it impossible to systematically study the behaviour of a single surface at various primes. When~the Picard rank is larger, prime numbers in a bigger range may be~used.

Concretely,~our sample consists of the resolutions of quartic surfaces having only
$A_1$~singularities.
We~chose 1600 quartic surfaces with 14 singularities. For~each of the surfaces,
we computed the upper bounds which were found at all the
primes~$p < 50$.
In~some cases, we continued the computations using larger primes up
to~$p = 103$.

It~turned out that good primes existed in every~example. We~could compute all the geometric Picard~ranks.
\end{ttt}

\begin{qu}
Do there exist good primes for all
$K3$~surfaces
over~$\bbQ$?
\end{qu}

\begin{lui}
\label{Luijk}
The~Picard~group of a
$K3$~surface
is isomorphic
to~$\bbZ^n$
where
$n$~may
range from
$1$
to~$20$.
An upper bound for the Picard~rank of a
$K3$~surface
may be computed as~follows. One~has the inequality
$$\rk \Pic (V_{\overline\bbQ}) \leq \rk \Pic (V_{\overline\bbF_p})$$
which is true for every smooth variety
$V$
over~$\bbQ$
and every
prime~$p$
of good~reduction~\cite[Example~20.3.6]{Fu}.

Further,~for a
$K3$~surface~$\calV$
over the finite
field~$\bbF_{\!p}$,
one has the first Chern~class~homomorphism
$$c_1\colon \Pic (\calV_{\overline\bbF_{\!p}}) \longrightarrow H^2_\et (\smash{\calV_{\overline\bbF_{\!p}}} \!, \bbQ_l (1))$$
into
$l$-adic~cohomology.
There~is a natural operation of the Frobenius on
$H^2_\et (\smash{\calV_{\overline\bbF_{\!p}}} \!, \bbQ_l (1))$.
All~eigenvalues are of absolute
value~$1$.
The~Frobenius operation on the Picard group is compatible with the operation on~cohomology.

Every~divisor is defined over a finite extension of the ground~field. Consequently,~on the subspace
$\smash{\Pic (\calV_{\overline\bbF_{\!p}}) \!\otimes_\bbZ\! \overline\bbQ_l \hookrightarrow H^2_\et (\calV_{\overline\bbF_{\!p}} \!, \bbQ_l (1))}$,
all eigenvalues are roots of~unity. These~correspond to eigenvalues of the Frobenius operation
on~$\smash{H^2_\et (\calV_{\overline\bbF_{\!p}} \!, \bbQ_l)}$
which are of the form
$p \zeta$
for
$\zeta$
a root of~unity. One~may therefore estimate the rank of the Picard
group~$\smash{\Pic (\calV_{\overline\bbF_q})}$
from above by counting how many eigenvalues are of this particular~form.

Doing~this for one prime, one obtains an upper bound
for~$\rk \Pic(V_{\overline{\bbF}_{\!p}})$
which is always~even. The~Tate conjecture asserts that this bound is actually~sharp. For~this reason, one tries to combine information from two~primes. The~assumption that the surface would have
Picard~rank~$2l$
over~$\overline{\bbQ}$
and~$\overline{\bbF}_{\!p}$
implied that the discriminants of both Picard groups,
$\Pic(V_{\overline{\bbQ}})$
and~$\Pic(V_{\overline{\bbF}_{\!p}})$,
were in the same square~class. Note~here that reduction
modulo~$p$
respects the intersection~product. When~combining information from two primes, it may happen that one finds the rank
bound~$2l$
twice, but the square classes of the discriminants are~different. Then,~these data are incompatible
with Picard
rank~$2l$
over~$\overline{\bbQ}$.
One~gets the rank
bound~$(2l - 1)$.
\end{lui}

\begin{rem}
There~are refinements of the method of van Luijk described in~\cite{EJ3} and~\cite{EJ5}. We~will not test these refinements~here.
\end{rem}

\begin{exa}
Let~$V$
be a
$K3$~surface
of
Picard~rank~$1$.
We~denote~by
$$V^n := \bigtimes_{i=1}^n V$$
the
\mbox{$n$-fold}
cartesian~product. Then,~the Picard rank
of~$V^n$
is equal
to~$n$.
Assuming~the Tate conjecture, one sees that the Picard rank of the reduction at an arbitrary prime is at
least~$2n$.

This~shows that there is no good prime
for~$V^n$.
Not~knowing the decomposition
of~$V^n$
into a direct product, we could not determine its Picard~rank.
\end{exa}

\paragraph{The analytic discriminant -- The Artin-Tate formula.}

For~the final step in~\ref{Luijk}, one needs to know the discriminant of the Picard~lattice.
One~possibility to compute this is to use the Artin-Tate~formula.

\begin{cono}[{\rm Artin-Tate}{}]
\label{AT}
Let\/
$V$
be a\/
$K3$~surface
over a finite
field\/~$\bbF_{\!q}$.
Denote~by\/
$\rho$
the rank and by\/
$\Delta$
the discriminant of the Picard group
of\/~$V\!\!$,
defined
over\/~$\bbF_{\!q}$.
Then,
$$|\Delta| = \frac{\lim\limits_{T \rightarrow q}\! \frac{\Phi(T)}{(T-q)^{\rho}}}{q^{21 -\rho} \#\!\Br(V)} \, .$$
Here,
$\Phi$
denotes the characteristic polynomial
of\/~$\Frob$
on\/~$\smash{H_\et^2(V_{\overline{\bbF}_{\!q}}\!, \bbQ_l)}$.
Finally,
$\Br(V)$
is the Brauer group
of\/~$V$.
\end{cono}

\begin{rems}
\begin{iii}
\item
The Artin-Tate conjecture is proven for most
$K3$~surfaces.
Most~notably, the Tate conjecture implies the Artin-Tate~conjecture~\cite[Theorem~6.1]{Mi1}. In~these cases,
$\#\!\Br(V)$
is a perfect~square.

On~its part, the Tate conjecture is proven for
$K3$~surfaces
under various additional~assumptions. For~example, it is true for elliptic
$K3$~surfaces~\cite{ASD}.
\item
In~such a case, the Artin-Tate formula allows to compute the square class of the discriminant of the Picard group over a finite~field. No~knowledge of explicit generators is~necessary.
\end{iii}
\end{rems}

\section{Singular quartics}

Singular quartic surfaces were extensively studied by the classical geometers
of the 19th~century, particularly by E.\,Kummer and~A.\,Cayley. For~example, the concept of a trope is due to this~period~\cite{Je}.

\begin{defi}
Let
$V \subset \bP^3$
be any quartic~surface. Then,~by a {\em trope\/}
on~$V$,
we mean a
plane~$E$
such that
$V \cap E$
is a double~conic. This~is equivalent to the condition that the equation
defining~$V$
becomes a perfect square
on~$E$.
\end{defi}

\begin{rem}
Tropes~lead to singular points on the surface
$V^\vee \subset (\bP^3)^\vee$
dual
to~$V$.
\end{rem}

\begin{lemo}[{\rm Kummer}{}]
A~quartic surface without singular curves may have at most 16 singular~points.
\eop
\end{lemo}

\begin{fac}
\label{ger}
Let\/~$V$
be a normal quartic~surface.

\begin{abc}
\item
Then,~not more than three singular points
on\/~$V$
may be~collinear.
\item
If~three singular points
on\/~$V$
are collinear then the line connecting them lies
on\/~$V$.
\end{abc}\smallskip

\noindent
{\bf Proof.}
{\em
b)
Otherwise,~this line would
meet~$V$
in each of the three singular points with multiplicity at least~two.\smallskip

\noindent
a)
Suppose~$k \geq 4$
singular points are~collinear. Then,~the line connecting them is contained
in~$V$.
Choose~a plane through this line not meeting any other~singularity.
The~intersection consists of the line and a possibly degenerate cubic~curve.
They~may not have more than three points in~common.
}
\eop
\end{fac}

\paragraph{A classical family.}
A~classification of the singular quartic surfaces with at least eight singularities of
type~$A_1$
was given by K.\,Rohn~\cite{Ro}. In~this article, we will deal with one of
the most important classical~families.

\begin{lemo}[{\rm Cayley, Rohn}{}]
\label{CR}
A family of quartics
in\/~$\bP^3$
such that the generic member has 14~singularities of
type\/~$A_1$
and no others is given by
\begin{eqnarray*}
\det
\left(
\begin{array}{cccc}
  0 &  l_1 &  l_2 &  l_3 \\
l_1 &    0 & l'_3 & l'_2 \\
l_2 & l'_3 &    0 & l'_1 \\
l_3 & l'_2 & l'_1 &    0
\end{array}
\right)
= 0 \, .
\end{eqnarray*}
Here,~$l_1, l_2, l_3, l'_1, l'_2, l'_3 \in \bbC[x,y,z,w]$
are linear~forms.
\end{lemo}

\begin{rem}
According~to K.\,Rohn~\cite{Ro}, every quartic surface with exactly 14 singular points, each being of
type~$A_1$,
is contained in this~family. Cf.~\cite[Ch.\,I, \S12]{Je}.
\end{rem}

\begin{rems}
\begin{iii}
\item
Evaluating the determinant, we find the explicit~equation
$$l_1^2 {l'_1}^2 + l_2^2 {l'_2}^2 + l_3^2 {l'_3}^2 - 2 l_1 l_2 l'_1 l'_2 - 2 l_1 l_3 l'_1 l'_3 - 2 l_2 l_3 l'_2 l'_3 = 0 \, .$$
\item
A~singular point is given
by~$l_1 = l_2 = l_3 = 0$.
As~the equation has the three independent symmetries
$l_i \leftrightarrow l'_i$,
there are eight singularities of this~type.

Two~further singular points are given
by~$l_1 = l'_1 = l_2l'_2 - l_3l'_3 = 0$.
As~the roles of the indices are interchangeable, there are a total of six singularities of this~form.
\item
Each of the six planes
$l_1 = 0, l_2 = 0, l_3 = 0, l'_1 = 0, l'_2 = 0, l'_3 = 0$
is a~trope. Generically,~these are the only tropes on such a~quartic.
Each~trope passes through six of the 14~singular~points.
\end{iii}
\end{rems}

\section{The desingularization}

\begin{lem}
\label{des}
Let\/~$\smash{\pi\colon \widetilde{V} \to V}$
be the desingularization of a normal quartic
surface\/~$V$
with only\/
$A_1$~singularities.
Then,~$\smash{\widetilde{V}}$
is a\/
$K3$~surface.\smallskip

\noindent
{\bf Proof.}
{\em
On~the smooth part
of~$V$,
the adjunction formula~\cite[Sec.~1.1, Example~3]{GH} may be applied as~usual. As,~for the canonical sheaf, one
has~$K_{\bP^3} = \calO(-4)$,
this shows that the invertible sheaf
$\Omega^2_{V^\reg}$
is~trivial.
Consequently,~$K_{\widetilde{V}}$
is given by a linear combination of the exceptional~curves.

However,~for an exceptional
curve~$E$,
we have
$E^2 = -2$.
Hence, according to the adjunction formula,
$K_{\widetilde{V}}E = 0$
which shows that
$K_{\widetilde{V}}$
is~trivial. The~classification of algebraic surfaces~\cite{Be} assures that
$\smash{\widetilde{V}}$
is either a
$K3$~surface
or an abelian~surface.

Further,~a standard application of the theorem on formal functions implies
$R^1\pi_* \calO_{\widetilde{V}} = 0$.
Hence,~$\chi_\alg(\widetilde{V}) = \chi_\alg(V) = 2$.
This~shows that
$\smash{\widetilde{V}}$
is actually a
$K3$~surface.
}
\eop
\end{lem}

\begin{rems}
\begin{iii}
\item
For~the assertion of the lemma, it is actually sufficient to assume that the singularities
of~$V$
are of types
$A$,
$D$,
or~$E$~\cite{Li}.
\item
In~general, the desingularization of a normal quartic surface is a
$K3$~surface,
a rational surface, a ruled surface over an elliptic curve, or a ruled surface over a curve of genus~three~\cite{IN}. The~latter possibility is caused by a quadruple~point. The~existence of a triple point implies that surface is~rational. It~is, however, also possible that there is a double point, not of type
$A$,
$D$,
or~$E$.
Then,~$\smash{\widetilde{V}}$
is rational or a ruled surface over an elliptic~curve.
\end{iii}
\end{rems}

\paragraph{Blowing up one singular~point.}
A generic line intersects a quartic
$V \subset \bP^3$
in precisely four~points. Assume~that
$P$~is
a double~point
on~$V$
which is not contained in a line lying
on~$V$.
Then,~the generic line
through~$P$
intersects
$V$
in two further~points. As~the lines
through~$P$ are parametrized by
$\bP^2$,
this leads to a double~cover
of~$\bP^2$
birational
to~$V$.

\begin{defi}
We~will call this scheme the {\em degree two model\/} corresponding
to~$V$.
\end{defi}

\begin{rems}
\begin{iii}
\item
It~is not hard to make the construction~explicit. For~this, suppose
that~$P = (0:0:0:1)$.
Then,~$V$
is given by an equation of the form
$Q(x,y,z)w^2 + K(x,y,z)w + F(x,y,z) = 0$
for a quadratic
form~$Q$,
a cubic
form~$K$,
and a quartic
form~$F$.
Multiplying~by~$Q$
and substituting
$W$
for~$Qw$
yields
$$W^2 + K(x,y,z)W + F(x,y,z)Q(x,y,z) = 0 \, .$$
The~ramification locus is the sextic curve given by
$4FQ - K^2 = 0$.
Actually,~when there are lines
through~$P$
lying
on~$V$,
this transformation works,~too.
\item
If~there is no line
on~$V$
containing~$P$
then the degree two model is simply the blow-up
of~$V$
in~$P$.
Indeed,~there is a morphism from the blow-up to the degree two model which is finite
and generically one-to-one. As~the degree two model is a normal scheme, Zariski's main theorem~applies.

In~general, the degree two model is the blow-up
of~$V$
in~$P$
with the lines containing
$P$
blown~down.
\item
Observe~that the conic
``$Q=0$''
is tangent to the ramification~sextic. Hence,~this conic splits in the double~cover. One~of the splits is actually the exceptional divisor produced by blowing up the singular~point.
\end{iii}
\end{rems}

\begin{rem}
When~we apply this construction to the particular singular quartics described above, the ramification sextic must have exactly 13 singular~points. According~to Pl\"ucker, such a highly singular degree-six curve is necessarily~reducible. It~is the union of three lines and a singular cubic or the union of two lines and two~conics.
\end{rem}

\section{Point~counting}

In~order to determine the eigenvalues of the Frobenius
on~$\smash{H_\et^2(\widetilde{V}_{\overline{\bbF}_{\!q}}\!, \bbQ_l)}$,
the usual method is to count the points
on~$V$
defined
over~$\bbF_{\!q}$
and extensions and to apply the Lefschetz trace~formula~\cite[Ch.\,VI, Theorem~12.3]{Mi}.

\begin{faco}[{\rm Elliptic fibration}{}]
\label{ef}
Let\/~$V \subset \bP^3$
be an irreducible quartic surface having at least two singular~points.
Then,~$\smash{\widetilde{V}}$
has an elliptic~fibration.\smallskip

\noindent
{\bf Proof.}
{\em
Intersect~$V$
with the pencil of hyperplanes through the two singular~points. This yields a fibration
of a surface birationally equivalent
to~$V$.
The~assumption implies that the generic fiber is an irreducible~curve. Depending~on
whether the line connecting the two singularities lies
on~$V$
or not, it is either a cubic curve or a quartic curve with at least two~singular~points.
In~both cases, Pl\"ucker's formulas show that its genus is at most~one. The~existence
of a fibration into curves of genus zero implied that
$\smash{\widetilde{V}}$
was~rational.
}
\eop
\end{faco}

\noindent
To~count the points on a singular quartic
surface~$V$
over~$\bbF_{\!q}$,
we have at least the following~possibilities.

\begin{algso}[{\rm Point counting}{}]

\begin{iii}
\item
Count~points~directly. This~means,
intersect~$V$
with a
\mbox{$2$-dimensional}
family of~lines. For~each line, determine the number points are on~it. This~last step means
to solve an equation of degree four
in~$\bbF_q$.
\item
Use~the elliptic~fibration. Enumerate~all fibers, defined
over~$\bbF_{\!q}$.
On~each fiber, count the number of~points.
\item
Compute~a degree~two model of the surface and count the points~there. This~means, on has to evaluate a sextic form
on~$\bP^2$
and to run an is-square routine in each~step.
\end{iii}
\end{algso}

\begin{rems}
\begin{abc}
\item
If the surface is defined
over~$\bbF_{\!p}$
then it suffices to count the points on a fundamental domain of the~Frobenius. This~leads to a
significant speed-up for all three~methods.
\item
In our examples, it turned out that the degree two model approach was the fastest~one. We~used it except for those surfaces where there were lines through each singularity defined
over~$\bbF_{\!p}$.
\end{abc}
\end{rems}

\section{Lower bounds for the Picard rank}

\begin{lem}
Let\/~$\smash{\pi\colon \widetilde{V} \to V}$
be the desingularization of a proper
surface\/~$V$
having only\/
$A_1$-singularities.

\begin{abc}
\item
Then,~the exceptional curves define a non-degenerate orthogonal system
in\/~$\smash{\Pic(\widetilde{V})}$.
\item
In~particular, the Picard~rank
of\/~$\smash{\widetilde{V}}$
is strictly bigger than the number of singularities
of\/~$V$.
\end{abc}\smallskip

\noindent
{\bf Proof.}
{\em
a)
The~exceptional curves have self-intersection
number~$(-2)$
each and do not meet each~other.\smallskip

\noindent
b)
For~$H$
the hyperplane section,
$\pi^* \calO_V(H)$
is orthogonal to the exceptional~curves.
}
\eop
\end{lem}

\begin{rem}
A~strategy to calculate the square~class of the discriminant
of~$\smash{\Pic(\widetilde{V})}$
is thus as~follows.
Consider~in~$\smash{\Pic(\widetilde{V})}$
the orthogonal complement
$P := \langle E_1, \ldots, E_n \rangle^\perp$.
Then,~$\smash{\disc \Pic(\widetilde{V}) \in 2^n (\disc P) \bbQ^* {}^2}$.
\end{rem}

\begin{ttt}
The~only method known to prove a non-trivial lower bound for the Picard~rank is to write down divisors~explicitly. We~always have the hyperplane~section. For~special quartics from the Cayley-Rohn~family, we observed two types of additional~divisors.

\begin{iii}
\item{\bf Lines.}
One~could search for lines on the surfaces by a Gr\"obner~base~calculation. However,~in our particular situation, every line connects at least two singular~points. We~will show this in Proposition~\ref{lines}~below.
\item{\bf Conics.}
There~is the special case that there exists a plane containing exactly four singularities no three of which are~collinear. Then,~the quartic curve on this plane splits into two~conics. The~same may happen when a plane through three singularities is tangent to the surface at another~point.
\end{iii}
\end{ttt}

\begin{ttt}
\label{intmat}
In~both these situations, one may directly calculate the corresponding intersection~matrices.

\begin{iii}
\item
Let~$k = 2, 3$
be the number of singularities connected by the
line~$l$.
Choose~a plane
through~$l$
such that the intersection curve splits into
$l$
and a smooth cubic~curve.
Then,~on~$\smash{\widetilde{V}}$,
we have two divisors
$L$
and~$C$
such that
$L^2 = -2$,
$C^2 = 0$,
and~$CE = 3 - k$.
For~$E_1, \ldots, E_k$
the exceptional divisors met
by~$L$,
$L' := L + \frac12 E_1 + \ldots + \frac12 E_k$
and
$C' := C + \frac12 E_1 + \ldots + \frac12 E_k$
are in the orthogonal complement of the exceptional~divisors (after tensoring
by~$\bbQ$).
Indeed,~this is an immediate consequence of Lemma~\ref{schnitt}, shown~below. We~find the intersection matrix
$-2 + k/2\;\; 3 - k/2 \choose \;\;\;3 - k/2\;\;\;\;\;\;\; k/2$
of
determinant~$2k - 9$.
\item
Here,~there are two conics meeting in four points
$k$~of
which are singular
on~$V$.
This~yields two divisors
$Q_1$
and~$Q_2$
on~$\smash{\widetilde{V}}$
such that
$Q_1^2 = Q_2^2 = -2$
and~$Q_1 Q_2 = 4 - k$.
In~a manner analogous to~i), we end up with the intersection matrix
$\!-2 + k/2\;\; \;\;\;4 - k/2 \choose \;\;4 - k/2\;\;\, -2 + k/2$
of
determinant~$2k - 12$.
\end{iii}
\end{ttt}

\begin{rem}
Observe~the following rules of thumb which apply as long as there are no multiplicities
$> \!\!1$~occurring.
If~$D$
meets exactly
$k$
singular points
then~$D'^2 = D^2 + \frac{k}2$.
If~$D_1 \neq D_2$
are irreducible curves having
$k$
singular and
$k'$
smooth points in common
then~$D_1 D_2 = k' + \frac{k}2$.
\end{rem}

\begin{lem}
\label{schnitt}
Let\/~$C$
be a curve on a
surface\/~$V$
having an\/
$A_1$-singularity
in\/~$P$.
Suppose\/~$P \in C$
and that\/
$P$~is
smooth
on\/~$C$.
Then,~on the
desingularization\/~$\smash{\pi\colon \widetilde{V} \to V}$,
the strict transform
of\/~$C$
meets the exceptional
curve\/~$\pi^{-1}(P)$
of order~one.\smallskip

\noindent
{\bf Proof.}
{\em
Indeed,~the model case for this situation is given by a conical
quadric~$V$
in~$\bP^3$
and a
line~$l$
on~$V$.
Then,~the desingularization is a Hirzebruch
surface~$\Sigma_2$
which is a ruled surface
over~$\bP^1$
with exactly one
\mbox{$(-2)$-curve}~$B$.
The~strict transform
of~$l$
is a
line~$F$
from the~ruling. It~is well known
that~$BF = 1$.
}
\eop
\end{lem}

\begin{propo}[{\rm Lines on special quartics}{}]
\label{lines}
Let\/~$V$
be a quartic surface with 14 singular~points. Then,~every line
on\/~$V$
contains two or three singular~points.\smallskip

\noindent
{\bf Proof.}
{\em
Let~$l$
be a line
on~$V$.
By~Fact~\ref{ger},
$l$
cannot contain more than three~singularities. Suppose~first that
$l$
is contained in one of the~tropes. Then,~this is a degenerate trope, the conic splitting into two~lines. As~a trope contains six singular points, there must be
three on each~line.

Otherwise,~$l$
meets each trope in a single~point. We~claim that these six points of intersection are all~singular. Then,~the assertion follows as a point can not be contained in more than three~tropes.

To~show the claim, assume that
$l$
would meet a trope in a smooth
point~$p$.
As~$l$
is supposed to be contained
in~$V$,
it is everywhere tangential
to~$V$.
But~for a point on a trope, the tangent~plane is the trope~itself.
Hence,~$l$
would be contained within the~trope. This~is a~contradiction.
}
\eop
\end{propo}

\section{Computations and numerical data}

\begin{ttt}
\label{nor}
Consider the Cayley-Rohn~family of determinantal quartics as described in Lemma~\ref{CR}.
Then,~over a Zariski open subset of the base, one may normalize to
$l_1 = x$,
$l_2 = y$,
$l_3 = z$,
and~$l'_1 = w$.
We~will write
$l'_2 = c_1 x + c_2 y + c_3 z + c_4 w$
and~$l'_3 = c_5 x + c_6 y + c_7 z + c_8 w$
for a coefficient
vector~$[c_1,\dots,c_8]$.
Over~a possibly smaller Zariski open subset, one has
$c_1, c_2, c_3 \neq 0$
in which case these coefficients may be normalized
to~$1$.
\end{ttt}

\paragraph{The computations carried out.}

We~chose a sample of 1600 singular quartics from the Cayley-Rohn~family. We~worked with the normal form as described in~\ref{nor}. The~coefficient vectors were produced by a random number~generator. The~coefficients themselves were integers in the
range~$-20, \ldots, 20$.
We~always
put~$c_1 = c_2 = c_3 = 1$.

For~each
surface~$V$
in the sample and each prime
$p < 50$
of good reduction, we counted the number of points
in~$V(\bbF_{\!p})$,
$V(\bbF_{\!p^2})$,
and~$V(\bbF_{\!p^3})$.
From~these data, we tried to compute the characteristic polynomial of the~Frobenius. For~the determination of the sign in the functional equation, we followed the strategy described in~\cite{EJ4}. We~used the explicitly known 15-dimensional sublattice of the Picard group generated by the hyperplane section and the exceptional curves in order to adapt the conditions to our~situation. In~430 cases, we had to compute, in
addition,~$\#V(\bbF_{\!p^4})$.
Here,~$p$
was up
to~$23$.
From~the characteristic polynomial, we read off the rank
of~$\Pic(V_{\overline\bbF_{\!p}})$
and, using the Artin-Tate~formula, computed its~discriminant.

\begin{rem}
In~the cases which remained with an unknown sign, we worked with the pair of possible characteristic~polynomials. This~means, we took the maximum of the predicted ranks as an upper~bound. In~the case that both upper bounds were equal
to~$16$,
we got a pair of possible square~classes for the~discriminants. Combining~information from different primes then meant to form the
intersection of these~sets.
To~give a typical example,
for~$p = 31$,
these sign problems occurred in 139 of 1299 cases with good~reduction. Other~primes led to similar~rates.
\end{rem}

\begin{rem}
According~to Fact~\ref{ef}, every surface in the sample is~elliptic. This~is enough to show that the Artin-Tate formula~\ref{AT} for the discriminant is~applicable.
\end{rem}

\paragraph{The average value for a prime.}

The~probability to obtain a good rank bound increases when the prime numbers~increase. Let~us visualize this by a~diagram.\vskip-1mm

\begin{figure}[H]
\centerline{
\includegraphics[height=8cm]{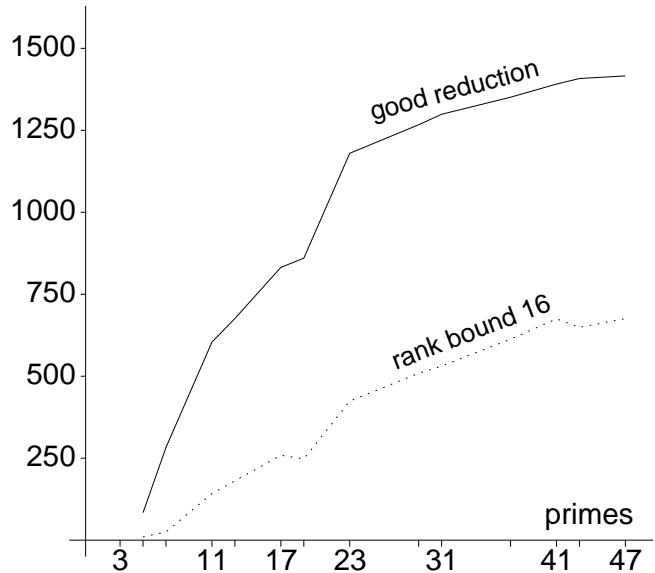}}\vskip-3mm
\caption{Number of surfaces with good reduction and rank bound 16 for $p < 50$.}
\end{figure}\vskip-10mm

\paragraph{The discriminants.}

We~computed the discriminant in all cases of Picard
rank~$16$.
In~4690 cases, we obtained a rank bound
of~$16$
and the determination of the sign in the functional equation was~possible. These~data led to 59 distinct square~classes for the~discriminant. The~most frequent square class was
$(-1)$
with 819~repetitions. The next one was
$(-2)$
having 608~repetitions. On~the other hand, each of the discriminants
$(-47)$,
$(-59)$,
$(-67)$,
$(-71)$,
$(-82)$,
$(-101)$,
$(-118)$,
$(-141)$,
$(-149)$,
and~$(-177)$
occurred only~once.

\paragraph{The ranks over~\boldmath$\bbQ$.}

On~each surface in the sample, we searched for additional~divisors. It~turned out that 1504 of the surfaces contained no line and no plane through four singular~points. For~these, we tried to prove that the Picard~rank
is~$15$.
For~the others, we tried to prove
Picard~rank~$16$.
The~statistics over the primes used is given by the table~below.

\begin{table}[H]
\scriptsize
\begin{multicols}{2}
\hfill
\columnsep3mm
\begin{tabular}{|c|c|c|}
\hline
prime  &  \#cases finished & \#cases left \\
\hline\hline
\phantom{0}11 & \phantom{00}2 &           1502 \\\hline
\phantom{0}13 & \phantom{0}15 &           1487 \\\hline
\phantom{0}17 & \phantom{0}36 &           1451 \\\hline
\phantom{0}19 & \phantom{0}57 &           1394 \\\hline
\phantom{0}23 &           151 &           1243 \\\hline
\phantom{0}29 &           181 &           1062 \\\hline
\phantom{0}31 &           219 & \phantom{0}843 \\\hline
\phantom{0}37 &           214 & \phantom{0}629 \\\hline
\phantom{0}41 &           173 & \phantom{0}456 \\\hline
\phantom{0}43 &           136 & \phantom{0}320 \\\hline
\phantom{0}47 &           118 & \phantom{0}202 \\\hline
\phantom{0}53 & \phantom{0}80 & \phantom{0}122 \\\hline
\phantom{0}59 & \phantom{0}44 & \phantom{00}78 \\\hline
\phantom{0}61 & \phantom{0}36 & \phantom{00}42 \\\hline
\phantom{0}67 & \phantom{0}20 & \phantom{00}22 \\\hline
\phantom{0}71 & \phantom{0}12 & \phantom{00}10 \\\hline
\phantom{0}73 & \phantom{00}6 & \phantom{000}4 \\\hline
\phantom{0}79 & \phantom{00}2 & \phantom{000}2 \\\hline
          103 & \phantom{00}1 & \phantom{000}1 \\
\hline
\end{tabular}
\vskip1.5mm

\centerline{~~~~~~~~~~~~~~~~~Rank 15 expected}
\columnbreak
\begin{tabular}{|c|c|c|}
\hline
prime  &  \#cases finished & \#cases left \\
\hline\hline
\phantom{0}5 &  \phantom{0}1  &            95 \\\hline
\phantom{0}7 &  \phantom{0}3  &            92 \\\hline
11           &  \phantom{0}3  &            89 \\\hline
13           &  \phantom{0}5  &            84 \\\hline
17           &  \phantom{0}2  &            82 \\\hline
19           &  \phantom{0}4  &            78 \\\hline
23           &            11  &            67 \\\hline
29           &  \phantom{0}7  &            60 \\\hline
31           &  \phantom{0}6  &            54 \\\hline
37           &  \phantom{0}8  &            46 \\\hline
41           &            12  &            34 \\\hline
43           &  \phantom{0}7  &            27 \\\hline
47           &  \phantom{0}6  &            21 \\\hline
53           &  \phantom{0}4  &            17 \\\hline
59           &  \phantom{0}1  &            16 \\\hline
61           &  \phantom{0}3  &            13 \\\hline
67           &  \phantom{0}3  &            10 \\\hline
73           &  \phantom{0}1  &  \phantom{0}9 \\\hline
79           &  \phantom{0}1  &  \phantom{0}8 \\\hline
83           &  \phantom{0}3  &  \phantom{0}5 \\\hline
97           &  \phantom{0}1  &  \phantom{0}4 \\
\hline
\end{tabular}\\\smallskip

\centerline{Rank 16 expected~~~~~~~~~~~~~~~~~}
\hfill
\end{multicols}\vskip-6mm
\caption{Progress of the upper~bounds}
\normalsize
\end{table}

\noindent
Observe~that there were a few cases where the data
for~$p < 50$
were not sufficient. For~these, we continued the point~count, in an extreme case up
to~$p = 103$.

\paragraph{Testing isomorphy.}

As~a byproduct of the computations, we proved that the surfaces in our sample are pairwise non-isomorphic. For~this, it was sufficient to show that, for each pair of surfaces, there existed a prime where both have good reduction, but the geometric Picard groups differ in rank or~discriminant. In~order to do this, we had to continue the point count in a few~cases. In~fact, the data
for~$p \leq 61$
contained enough~information.

\paragraph{The five examples left.}

\begin{exa}
Let~$S_1$
be the surface given by the coefficient~vector
$[ 1, 1, 1, -7, 16, 6, -9, 12 ]$.
Here,~there is a plane through three singularities which is tangent
to~$S_1$
at a fourth~point. The~intersection curve splits into two~conics. The~Picard rank is thus at
least~$16$.
On~the other hand, we found the rank
bound~$16$
for~$p =  61$,
$71$,
$83$,
and~$101$.
\end{exa}

\begin{exas}
Let~$S_2$,
$S_3$,
and~$S_4$
be the surfaces given by the vectors
$[ 1, 1, 1, -1, -16, 7, 10, -10 ]$,
$[ 1, 1, 1, 3, -16, 2, 4, 15 ]$,
and
$[ 1, 1, 1, -1, 13, -11, 1, 15 ]$,
respectively. For~each surface, we found
rank~$18$
at several primes with various~discriminants. Hence,~in each case, there was an upper bound
of~$17$
for the Picard~rank.

\begin{iii}
\item
On~$S_2$,
we found a
plane~$E$
through four singular
points~$P_1, \ldots, P_4$.
On~$E$,
the quartic splits into two
conics~$Q_1, Q_2$.
Further,~there are two
lines~$L_1, L_2$
through
$P_1$
and~$P_2$
on~$V$
which meet in a smooth~point. Actually,
$L_1$
and~$L_2$
form a degenerate~trope. Arguing~as in~\ref{intmat}, we find the intersection~matrix
$$\left(
\begin{array}{rrrr}
      0 & \phantom{-}2 &  \frac12 &  \frac12 \\
      2 &            0 &  \frac12 &  \frac12 \\
\frac12 &      \frac12 & -\frac12 &        1 \\
\frac12 &      \frac12 &        1 & -\frac12
\end{array}
\right)$$
of rank~three, This confirms Picard
rank~$17$.
\item
On~$S_3$,
the situation is analogous to that
on~$S_2$.
The~only difference is that the
plane~$E$
meets three singular points and is tangent
to~$V$
at a fourth~point.
$L_1$~and~$L_2$
meet~$E$
in singular~points. We~find the intersection~matrix
$$\left(
\begin{array}{rrrr}
-\frac12 &  \frac52 &  \frac12 &  \frac12 \\
 \frac52 & -\frac12 &  \frac12 &  \frac12 \\
 \frac12 &  \frac12 & -\frac12 &        1 \\
 \frac12 &  \frac12 &        1 & -\frac12
\end{array}
\right)$$
of rank~three. Again,~this confirms Picard
rank~$17$.
\item
Here,~we have a
plane~$E$
through five singular~points.
On~$E$,
the quartic splits into a
conic~$Q$
and two
lines~$L_1, L_2$.
There~are two further lines
$L_3, L_4$
through three~singularities.
$L_1$~and
$L_3$
meet in a smooth~point. Together,~they form a degenerate~trope. The~same is true for
$L_2$~and
$L_4$.
Finally,~$L_3$
and~$L_4$
have a singular point in~common. We~find the intersection~matrix
$$\left(
\begin{array}{rrrrr}
0 &        1 &        1 &        0 &        0 \\
1 & -\frac12 &  \frac12 &        1 &        0 \\
1 &  \frac12 & -\frac12 &        0 &        1 \\
0 &        1 &        0 & -\frac12 &  \frac12 \\
0 &        0 &        1 &  \frac12 & -\frac12
\end{array}
\right)$$
of rank~three. Again,~this confirms Picard
rank~$17$.
\end{iii}
\end{exas}

\begin{exa}
Let~$S_5$
be the surface given by the coefficient~vector
$[ 1, 1, 1, -1, -13, 0, 11, -11 ]$.
We got a rank bound
of~$18$
for~$p = 23$,
$31$,
$61$,
$79$,
$89$,
$97$,
and~$101$.

On~the other hand, we found quite a number of particular divisors on this~surface. There~are six planes through exactly four singular~points. On~each of these planes, the quartic splits into two~conics. The~combinatorial structure is rather~interesting. In~a somewhat arbitrary numbering, the table below describes which plane meets which singular~points.

\begin{table}[H]
\scriptsize
\begin{center}
\begin{tabular}{|c||c|c|c|c|}
\hline
$E_1$ & $P_1$ & $P_2$ & $P_9$    & $P_{11}$ \\\hline
$E_2$ & $P_3$ & $P_4$ & $P_{13}$ & $P_{14}$ \\\hline
$E_3$ & $P_3$ & $P_7$ & $P_{12}$ & $P_{14}$ \\\hline
$E_4$ & $P_4$ & $P_8$ & $P_{10}$ & $P_{13}$ \\\hline
$E_5$ & $P_5$ & $P_6$ & $P_9$    & $P_{11}$ \\\hline
$E_6$ & $P_7$ & $P_8$ & $P_{10}$ & $P_{12}$ \\\hline
\end{tabular}
\end{center}\vskip-5mm
\caption{Planes through four singular points}
\normalsize
\end{table}

\noindent
Further,~there are the two lines
$L_1$
through
$P_9, P_{10}$
and~$P_{14}$
and
$L_2$
through
$P_{11}, P_{12}$
and~$P_{13}$.
The~lines
$L_1$
and~$L_2$
have a smooth point in~common. They~form a degenerate~trope. The~intersection matrix of
$L_1$,
the two conics
in~$E_2$,
and one of the conics
in~$E_3$
alone~is
$$\left(
\begin{array}{rrrr}
-\frac12 &      \frac12 &      \frac12 &      \frac12 \\
 \frac12 & \phantom{-}0 & \phantom{-}2 & \phantom{-}1 \\
 \frac12 &            2 &            0 &            1 \\
 \frac12 &            1 &            1 &            0
\end{array}
\right)$$
of rank~four. This~proves that the Picard rank
is~$18$.
\end{exa}

\begin{summ}
We~considered the resolutions of 1600 randomly chosen Cayley-Rohn quartics with exactly 14 singularities of
type~$A_1$.
The~corresponding
$K3$~surfaces
were mutually non-isomorphic. It~turned out that all the Picard ranks could be~determined. However,~at several examples rather large primes up
to~$p = 103$
had to be~considered. We~found Picard rank~fifteen
$1503$~times
and Picard rank sixteen
$93$~times.
Further,~there were three surfaces of Picard rank~seventeen and one surface of Picard rank~eighteen in the~sample.
\end{summ}

\end{document}